# [1]On Mathematical Symbols in China

By    Fang Li, Yong Zhang

(Department of mathematics, Zhejiang University, Hangzhou)

When studying the history of mathematical symbols, one finds that the development of mathematical symbols in China is a significant piece of Chinese history; however, between the beginning of mathematics and modern day mathematics in China, there exists a long blank period. Let us focus on the development of Chinese mathematical symbols, and find out the significance of their origin, evolution, rise and fall within Chinese mathematics.

## 1. The origin of the mathematical symbols

The symbols for numerals are the earliest mathematical symbols. Ancient civilizations in Babylonia, Egypt and Rome, which each had distinct writing systems, independently developed mathematical symbols. Naturally, China did not fall behind. In the 16th century BC, symbols for numerals, called Shang Oracle Numerals because they were used in the Oracle Bone Script, appeared in China as a set of thirteen numbers, seen below. (The original can be found in reference [2]).

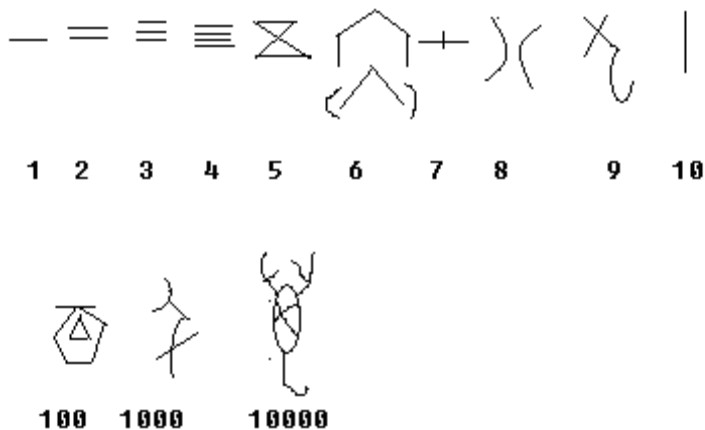

**Figure 1**

Oracle Bone Script is a branch of hieroglyphics. In Figure 1, it is obvious that 1 to 4 are the hieroglyphs for numbers, while the remaining numbers are phonetic loans from the names of animals and plants, the most

---

obvious of which are the symbols for "hundred", named after the Chinese word for pinecone and "ten thousand", named after the scorpion (some argue that the word is a result of the scorpion's large number of feet). Some experts maintain that the word "hundred"（百）which was also the symbol for "speak continuously" in ancient China, was derived from the word "white"（白）which symbolizes the human's head in hieroglyphics. Similarly, the Chinese word "ten thousand"（万）in Oracle Bone Script was derived from the symbol for scorpion, possibly because it is a creature found throughout rocks "in the thousands"。 In Oracle Script, the multiples of ten, hundred, thousand and ten thousand can be denoted in two ways: one is called co-digital or co-text, which combines two single figures; another one is called analysis-word or a sub-word, which uses two separate symbols to represent a single meaning（see [5]）. In addition to the Oracle Script, there was a script created in the Western Zhou Dynasty (1046 BC – 771 BC), called "Vase Character" or "Jin Character", referring to the characters that were engraved in bronze utensils like vases.

## 2. From Counting Rod to Abacus

It has been a long time since numbers were recorded by nicks and objects, a method that led to the invention of counting symbols. In ancient China, the main counting and calculation tool was the counting rod, the use of which was said to be "Single digits use vertical rods, tens digits take horizontal rods, hundreds re-use vertical rods, and then thousands re-use horizontal, ten thousands re-use vertical rods...... so from right to left, vertical rods and horizontal rods appear by turns, and so on,; you can express with the counting rod extremely large natural numbers". From the end of Shang Dynasty to the end of the Western Zhou Dynasty, China gradually developed the modern counting rod. It has been over two thousand years since the counting rod was first invented in China. The counting rod is actually a decimal counting method, which proves that China held a leading place among the ancient civilizations in the use of the decimal counting technique. The principle is to "know bits first", using vertical and horizontal rods. To signify a zero, the position is left blank. There are two ways to denote numbers in the counting rod system, the vertical and the horizontal (as shown):

```
vertical
type   :  |  ||  |||  ||||  |||||  T  ⊤  ⊤⊤  ⊤⊤⊤
horizontal
type   :  —  =  ≡  ≣  ≣≣  ⊥  ⊥  ⊥⊥  ⊥⊥⊥

number :  1   2   3   4    5    6   7   8    9
```

Later this way of denoting numbers was named the "counting rod digitals". Since then, number scripts have evolved continuously; among them the most outstanding is the "Southern Song Digital", which was created by Qin Jiushao (Courtesy Name Daogu), a mathematician who lived in the Southern Song dynasty. In his book —*Mathematic Treatise in Nine Sections*—, the "Southern Song Digital" has become a specific counting symbol throughout the whole book. Like the counting rod, the "Southern Song Digital" contained two forms, upright and horizontal:

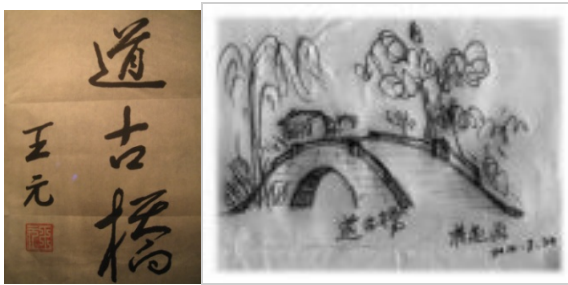

**On the Xixi Road (which is between Hangda Road and Huanglong Road, closer to Hangda), there is a stone bridge named Daogu Bridge and it was first built in the Southern Song Dynasty of Jiaxi years (1237-1241) and originally named the Xixi Bridge. It was recorded in Xian Chun's *Linan Load* of the Southern Song Dynasty, "The Xixi Bridge is on the east of the Government's Trial Court and was built in the Song Dynasty of Jiaxi years by Daogu". The builder Daogu was indeed the mathematician Qin Jiushao. [1]**

```
vertical
type    |  ||  |||  ||||  X  |||ōȯ  T  ⊤  ⊤⊤  ⊤⊤⊤X  O

horizontal
type    —  =  ≡  ≣  X  ||||ōȯ  ⊥  ⊥⊥  ⊥⊥  ⊥⊥⊥X  O

number  1  2  3  4  5  6  7  8  9  0
```

Among them, ||||X and ≡X are used to replace |||| and ≡ for representing "4", ||||ōȯ and ≡ōȯ are used to replace |||| and ≡ for representing "5", and ⊤⊤⊤X and ⊥X are

used to replace 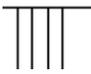 and 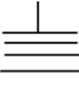 for representing "9". Addition, subtraction, multiplication, division, extracting square root, extracting cubic root, and finding a root of an algebraic equation of certain degrees can be accomplished by manipulating counting rods on a board. Some of ancient China's mathematic achievements were completed in the counting rod era such as Zu Chongzhi's calculation of $\pi$, equation solutions and extraction of a root in Qin Jiushao's—'*The Nine Chapters on the Mathematical Art*'—. However, there were limitations to the counting rod. For example, general polynomial equations with five or more variables could not be solved. Therefore from the Three Kingdom period to the late Tang Dynasty, there was a transition from the ancient counting rod to the abacus. During the Ming Dynasty, the counting rod was gradually replaced by the abacus on the stage of Chinese mathematics.

Compared with the counting rods, the abacus overcame the disadvantages of the upright and horizontal numeration system and inconvenience of placing the rods while having the advantage of a decimal place system like the counting rod system. However, when the abacus was first invented, its calculation beads were not stringed together, which made it inconvenient to carry and restricted its wide adoption. In the Northern Song Dynasty, the abacus with stringed beads was famous throughout the country; but the algorithms and formulas were not developed until the Yuan Dynasty.

In the 1910s, the India – Arabic numerals took place of the abacus and are still in use today. In fact, as early as the beginning of the 13th century, the India-Arabic numerals were introduced to China, but were not widely adopted.

## 3. "Tianyuan system" and "Siyuan system"

The Song and Yuan Dynasties were the most prosperous times in ancient China. During these periods, one of the most profound trends was the attempt to use algebra symbols, leading to the invention of the "Tianyuan system" and "Siyuan system" (see [2]). They were methods of denoting variables with special symbols to solve equations. The "Tianyuan system" first appeared in *Mathematic Treatise in Nine Chapters* written by the Southern Song Dynasty mathematician, Qin jiushao (Courtesy name: Daogu). It first appeared systematically in *Sea Mirror of Circle Measurement* written by Li Ye in the Jin and Yuan Dynasties and it was the first time mathematical symbols were introduced in Chinese mathematical history, earlier than similar algebra in Europe by a few centuries. The distinguished mathematician in the Yuan Dynasty, Zhu Shijie, further illustrated the advantage of the "Tianyuan system". To write out an equation with "Tianyuan system", one should first "suppose Tianyuan is something", that is, "Tianyuan" is the variable equivalent to "suppose x is something" today; then one should give an equality of two polynomials about "Tianyuan" and moreover subtract two polynomials from each other to get an equation of polynomial of higher degree; finally, one gets the positive root of the equation through radication. Zhu Shijie discussed the "Siyuan method" systematically in *Precious Mirror of the Four Elements* (see [4]). In this book, Zhu Shijie proposed that in addition to "Tianyuan", "Diyuan(地元)", "Renyuan(人元)" and "Wuyuan(物元)" could also be variables if necessary and then gave the method of placing the counting rod of variables and constants. Also, he explained with examples how to eliminate the variables in multivariate set of equations and get a high degree polynomial equation with only one variant. Some argued that the "Siyuan system" is the most significant achievement and glorious chapter in ancient Chinese mathematics, in particular, in the area of research of equations, and the most outstanding in the realm of medieval mathematics all over the world. Unfortunately, the "Tianyuan system" and "Siyuan system" are regarded as imperfect symbol systems, because they still could not denote sets of equations with five variables or more.

The Siyuan system and Tianyuan system are two of the major achievements of mathematics in China, and the reason they did not develop further is because there was no unified system of symbols. Nevertheless, they still embody the preliminary exploration and research on mathematical symbols by Chinese people, because such narrative language itself is a kind of symbolic expression.

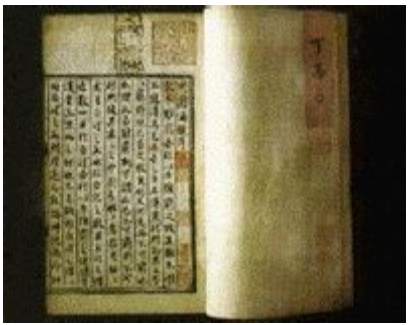

*Sea Mirror of Circle Measurement*

## 4. Regret of Mathematical symbols in China

Chinese traditional culture attached more importance to literature than science, which resulted in a low position of mathematicians in society. Also, because "one could not specialize in mathematics", mathematicians worked by themselves. Naturally, it was hard to form a mathematic community in a real sense and a uniform mathematical symbol system, which led to chances for the improvement of Chinese mathematics and mathematic symbols to be missed again and again.

Ancient China was one of the first civilized countries in the world to recognize and apply decimals. The conception of decimals in China can be dated to the 3rd or 4th century A.D. In *Master Sun's Mathematical Manual*, it is said that the calculation of length is started with hu (忽), ten hu is one si(丝), ten si is one hao(毫), ten hao is one li(厘) and ten li is one fen(分), with fen as the boundary of integers and decimals. That could be regarded as the origin of decimal name and numeration. It is proposed in *The Nine Chapters on the Mathematical Art* that if a number cannot be radiated finitely, it is irradiated and is numbered with "mian(面)"; this expressed the earliest recognition of irrational numbers in Chinese.

During the Han and Wei Dynasties and the Three Kingdoms period, the great mathematician Liu Hui appeared and proposed the concept of the decimal. When he was making commentaries for *The Nine Chapters on the Mathematical Art*, he concluded that if a number has a remainder after radicating, its radication may be represented in the form of decimal fraction (i.e. decimal). As he explained in *The Nine Chapters on the Mathematical Art*: "When one cannot extract a square root, then the root is expressed as a decimal. The figures following the decimal point are represented as a tenth if there is one and as a hundredth if there are two, and so on". Unfortunately, this numerical method did not develop further into a systematical theorem and so-called "decimal point symbol". Nevertheless, in ancient China the decimal and the decimal point were denoted with characters which could be shown in the document [3]: "cheng fen" denoted decimals. For example, 0.174 could be read as "one cheng seven fen four centimeters"; today, "cun(寸)" denotes the decimal point.

In addition to decimals, China was one of the first countries to discover decimal fractions. In ancient times, there were usually two numerical symbols for fractions: in characters, fractions were denoted as "tenths of"; with the counting rod method, the dividend was denoted as "shi"(实) at a middle place and, with the divisor denoted as "fa"(法) on the bottom and the quotient on the top.

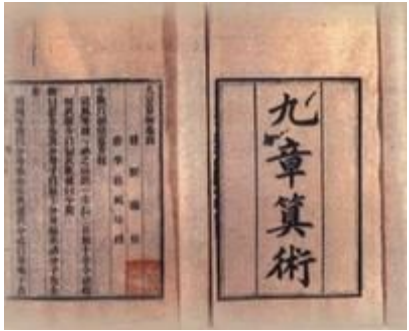 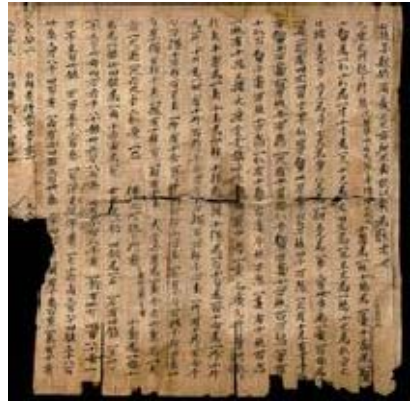

Master Sun's Mathematical Manual

China was also the first country in the world to discover and understand negative numbers. Liu Hui once said: "There are two opposite kinds of counting rods for gains and losses which are called positive and negative respectively, distinguished by upright and slanting rods respectively". In ancient China, there were many methods to denote negative numbers. "Red counting rods are positive and black ones are negative" means that ancient people used red for positive and black for negative, in addition to upright rods for positive and diagonal ones for negative.

As with their invention of fractions and decimals, the ancient Chinese stepped into the boundary of invention of negative numeration but fell back and again missed the chance to create new symbols for humankind. There are many reasons for this phenomenon. First of all, we think that as an ingrained numeration method in ancient times, the counting rod was widely used and only the algorithm, calculation data and calculation results were recorded in books without a reasoning process or calculation symbols. Additionally, in the Chinese feudal society, mathematics was mainly used to measure the land, launch water-conservation projects, allocate the labor force and calculate the collection of food, uses which are called "mathematics for management" nowadays. For example, in ancient China the method to calculate figure areas was produced directly from the measurements of field. In the chapter, Field Measurement, of *The Nine Chapters on the Mathematical Arts*, the calculation method for almost all figures are introduced, where the figures are named according to the shape of fields, like the square field (方田), the circle field (圆田), the arc field (弧田), the triangle field (圭田) and the trapezoidal field (邪田).

## 5. Mathematical symbols "bloom in doom" in China

Mathematics has existed as an art since ancient China. It was listed as one of the Six Arts in the Zhou Dynasty and regarded as a compulsory course for the royal children. During the Tang Dynasty, *The Nine Chapters on the Mathematical Arts* and other nine-section mathematical books were designed to be compulsory learning for cultural study, and mathematics was listed as one of the subjects examined in the general examination. In the Song and Yuan Dynasties, ancient mathematics had reached its peak, as mentioned above. However, after the Ming Dynasty, ancient China entered into the late imperial period. The government carried out authoritarian rule and the "eight-legged essay" (A style of essays used in the civil examination). The development of mathematics gradually declined from then on. In the late Ming Dynasty and the early Qing Dynasty, Western missionaries came to China and introduced Western scientific knowledge into China. Mathematical symbols, as a mathematical language, were also brought in with the introduction of mathematics. However, due to the closure policy in the Qing Dynasty, the introduction of Western knowledge was stopped for some time and was not reorganized until the end of the Opium War. After this, lots of Western works were translated into Chinese, and modern mathematics entered into Chinese mathematician' horizon, along with modern mathematic symbols.

However, in the Qing Dynasty, with literature greatly valued, the mathematicians, although they accepted the new ideas of the West, adhered to cultural tradition and generally believed that denoting with characters was superior to denoting with abstract symbols. As stated in [3], this phenomenon was especially obvious when the Western works were translated. They used characters to substitute for Western symbols, and even took pains to create new characters and a strange symbol system, which made the algebra too hard to recognize and understand. For example, in the book *Elements of Analytical Geometry, and of the Differential and Integral Calculus* translated together by Li Shanlan, the famous Chinese mathematician and translator, and Alexander Wylie, an English missionary, the Arabic digitals "1, 2, 3..." were replaced by Chinese characters; the English letters "a, b, c..." were replaced by Chinese "first, second, third..." and other Chinese characters; a capital letter was denoted by its lowercase letter with a mouth(口) as a component; "⊥" represented "+"; "⊤" represented "−"; the characters for "Differential" and "Integral" represented respectively the differential symbol "d" and the integral symbol "$\int$" and the sum symbol was replaced by the Chinese character "口昂". In *Physic Calculation*, he replaced the 28 Greek letters with 28 star names in Chinese, translated "$\pi$" into "Zhou" and wrote fractions with the numerator on the bottom and the denominator on the top. Although, from our point of view, these ridiculous symbols are obtrusive, they were viewed as model symbols at that time. All of these changes set up barriers to the integration of Chinese mathematics into world mathematics; However, they still played a positive role to some extent and paved the way to international mathematics since they broke the tradition of only using Chinese characters in mathematics. During the Xinhai Revolution, this unattractive symbol system finally had "its bones buried", and Chinese mathematics accomplished its complete symbolization and started to step into the family of world mathematics.

Ancient Europeans struggled with fractions because of the use of tedious Roman numerals, so ancient Europeans were weakest in algebra. However, when Indian-Arabic numerals were introduced into Europe, four operations in Europe were simplified. In addition to a series of simple symbol application. as a result, algebra developed rapidly. Similarly, China's ancient mathematics failed to form a perfect symbol system and were stagnant because of the limits of the counting rod, resulting in not only a lack of the axiomatic system, but also the absence of symbolic algebra, analytic geometry, calculus and modern variable mathematics in China. Naturally, to move forward it was necessary to break barriers and to get in touch with international standards and to forward to symbolic mathematics. This shows that the extension and development of any prosperous culture cannot be confined to the boundaries of that nation, and it is hard to create a complete and rich symbol system independently and separately from the influence of other nations and the outside world.

As the saying goes, "April showers bring May flowers". Chinese mathematics, which has a glorious tradition and has gone through setbacks and stagnation, is now integrated with international mathematics development in both its symbol system and research themes, and it has overcome its historical weaknesses. In our opinion, there is no doubt that the strong modern mathematical symbol system provided a powerful boost to original mathematical ideas in China.

Acknowledgements:

The authors thank Jie Dong and Xianjun Liu for their important proofread.

This project is supported by the National Natural Science Foundation of China (No. 11271318, No. 11171296 and No. J1210038) and the Specialized Research Fund for the Doctoral Program of Higher Education of China (No. 20110101110010) and the Zhejiang Provincial Natural Science Foundation of China (No.LZ13A010001)